\documentclass[twoside,leqno,12pt]{article}
\usepackage{ifthen,amsmath,amssymb,latexsym,graphicx,epsfig}
\numberwithin{equation}{section} \numberwithin{lemma}{section}
\numberwithin{theorem}{section} \numberwithin{corollary}{section}
\numberwithin{proposition}{section}
\numberwithin{definition}{section} \numberwithin{example}{section}
\numberwithin{remark}{section}
\def\TITLE       {ABOUT 
THE DIOPHANTINE EQUATION $x^{4}-q^{4}=py^{r}$}
\par
\def\LASTNAMEI   {\textbf{Savin}}
\def\FIRSTNAMEI  {\textbf{Diana}}

\def\ABSTRACT {
\begin{center}
\textbf{ABSTRACT}:
\end{center}
In this paper, we prove a theorem about the integer solutions to the Diophantine equation $x^{4}-q^{4}=py^{r}$, extending previous work of K.Gy\H ory, and F.Luca and A.Togbe, and of the author.\ \

}
\def\CLASSIFICATION{MSC (2000): 11D41}

\def\KEYWORDS{KEYWORDS: Diophantine equations; Kummer fields; cyclotomic fields}

\begin{document}

\begin{center}
\par\TITLE
\par\FIRSTNAMEI \ \LASTNAMEI
\par\medskip

\end{center}
\par\
\par\ABSTRACT{}
\par\
\par\CLASSIFICATION{}
\par\KEYWORDS{}
\par\


 \section{Introduction}
 \label{Introduction}
 Many Diophantine equations have been studied connected with the one in the title. For example,K\' alm\' an Gy\H ory studied (in [4],[5]) the Diophantine equation  $x^{p}+y^{p}=cz^{p}.$\\
 B. Powell and Henri Darmon studied the Diophantine equation $x^{4}-y^{4}=z^{p}.$\\
\indent In [10] B. Powell proved that this equation has no integer solutions with $p$ does not divide $xyz.$ In [3] H. Darmon obtained (using elliptic curves) the following result:\\
 \textit{Let} $p\geq3$ \textit{be a prime.Then:}\\
 i) \textit{the Diophantine equation} $x^{4}-y^{4}=z^{p}$ \textit{has no nontrivial solutions if} $p\equiv1$$\pmod4.$\\
 i) \textit{the Diophantine equation} $x^{4}-y^{4}=z^{p}$ \textit{has no nontrivial solutions with} $z$ \textit{even.}\\
 \smallskip\\
\indent In some previous papers ([12],[13],[14]) we considered some Diophantine equations of the form $x^{4}-q^{4}=py^{r},$ with $r\in\left\{3,5,7\right\}$, where $p,q$ are distict prime natural numbers satisfying some conditions.  Florian Luca and Alain Togbe have recently studied the equation from the title in the case $r=3.$ In [9], they showed that
the Diophantine equation $x^{4}-q^{4}=py^{3}$
has no integer solutions $(x,y,p,q)$ with $\gcd(x,y)=1,$ $xy\neq0,$ and$p$ and primes.\\
\indent Here we try to generalize the results from the papers [12],[13], taking $r$ a prime natural number different from $p$ and $q$, all of them satisfying conditions which will be given.\\
\indent The main result in this paper is:\\
 \smallskip\\
 \textbf{Main Theorem.}\textit{Let} $p,q,r$ \textit{be distinct prime numbers satisfying the conditions:\\
 $q\neq2,$ $p\equiv3$$\pmod4$, $p\equiv1$$\pmod r$, $r\equiv \pm 3$$\pmod8$, $\overline{p}$ is a generator of the group $\left(U(\mathbb{Z}_{q^{r-1}}),\cdot\right)$, $\overline{q}$ is a generator of the group $\left(\mathbb{Z}_{r}^{*},\cdot\right)$, $2$ is an r-power residue mod $q$.
  Then, any solution in coprime integers} $(x,y)$ \textit{to the equation}\\
   $$x^{4}-q^{4}=py^{r}$$
 \textit{satisfies the property that} $p$ \textit{must divide} $y$.
 \section{Preliminaries}
 \label{Preliminaries}
 \indent The proofs involve techniques based on the theory of Kummer fields and cyclotomic fields. For convenience sake which we recall in this section those properties of ideals in integer rings of such fields which we will be using in our proofs. 
 \smallskip\\
For a prime number $p$ and for $\zeta$ a primitive root of order $l$ of unity (where $\gcd(p,l)=1$), we have the following proposition:\\
\smallskip\\
\textbf{Proposition 2.1.}($\left[2\right]$).\textit{Let} $l\geq3$ \textit{be a positive integer}\textit{and} $\zeta$ \textit{a primitive} $l$-\textit{th root of unity}. \textit{Let} $\mathbb{Z}\left[\zeta\right]$ \textit{be the ring of integers of the cyclotomic field} $\mathbb{Q}\left(\zeta\right)$. \textit{If} $p$ \textit{is a prime natural number}, $l$ \textit{is not divisible by} $p$, \textit{and} $f$ \textit{is the smallest positive integer such that} $p^{f}$$\equiv$$1$ (\textit{mod}\:$l$), \textit{then we have}\\
$$p\mathbb{Z}\left[\zeta\right]=P_{1}P_{2}...P_{r},$$\\
 \textit{where} $r=\frac{\varphi\left(l\right)}{f}$, $P_{j},$ $j=1,2,...,r$
 \textit{are different prime ideals in the ring} $\mathbb{Z}\left[\zeta\right]$.\\
\smallskip\\
For the ring of integers in the Kummer field $\mathbb{Q}\left(\sqrt[l]{a};\zeta\right)$, where $a$$\in$$\mathbb{Z}$, the ideal $PA$, with $P$$\in$ Spec($\mathbb{Z}\left[\zeta\right]$), is totally characterized by the $l$ power-character
$\left\{\frac{a}{P}\right\}$, as in the following theorem. 
\smallskip\\
 \textbf{Theorem 2.2.}($\left[ 7\right],\left[ 8\right]$).\textit{Let} $l$ \textit{be a prime number,} $\zeta$ \textit{a primitive} $l$-\textit{th root of unity and} $A$ \textit{the ring of integers of the Kummer field} $\mathbb{Q}\left(\sqrt[l]{a};\zeta\right),$ \textit{where} $a\in$$\mathbb{Z}$, \textit{and} $P$ \textit{be a prime ideal in the ring} $\mathbb{Z}\left[\zeta\right]$. \textit{Then the following statments hold:}\\
  i) \textit{The ideal}\ $PA$ \textit{is equal with the l-power of a prime ideal in the ring} $A$, \textit{if} $\left\{\frac{a}{P}\right\}=0.$\\
  ii)\textit{The ideal}\ $PA$ \textit{decomposes in l different prime ideals in the ring} $A$, \textit{if} $\left\{\frac{a}{P}\right\}=1.$\\
  iii)\textit{The ideal}\ $PA$ \textit{is a prime ideal in the ring} $A$, \textit{if} $\left\{\frac{a}{P}\right\}$ 
  \textit{ equals a root of order l of unity, different from} $1$.
\smallskip\\
Now we recall properties of some Galois extensions of Kummer or cyclotomic types.  
\smallskip\\
  \textbf{Proposition 2.3.}($\left[ 11\right]$).\textit{Let A be the ring of integers of the Kummer field} $ \mathbb{Q}\left(\sqrt[l]{a};\zeta\right),$ \textit{where} $a$ \textit{is a positive integer and let} $\zeta$ \textit{be a primitive} $l$-\textit{th root of unity}. \textit{Let} $G$ \textit{be the Galois group of the Kummer field} $ \mathbb{Q}\left(\sqrt[l]{a};\zeta\right)$ \textit{over} $\mathbb{Q}\left(\zeta\right)$. \textit{Then, for any} $\sigma$$\in$$G$ \textit{and for any} $P$$\in$Spec$\left(A\right)$, \textit{we have} $\sigma\left(P\right)$$\in$Spec$\left(A\right)$.\\
\smallskip\\
\textbf{Proposition 2.4.}($\left[ 11\right]$).\textit{Let be given an extension of fields}  $\mathbb{Q}$$\subset$ $\textbf{K}$=$\mathbb{Q}\left(\zeta,\sqrt[l]{a}\right)$,\textit{where} $\zeta$  \textit{is a primitive} $n$-\textit{th root of unity}. \textit{Then the extension} $\mathbb{Q}$$\subset$$\textbf{K}$ \textit{is a Galois extension, the Galois group} $G(\textbf{K}/\mathbb{Q})$ \textit{is solvable and the Galois group} $G(\textbf{K}/\mathbb{Q}\left(\zeta\right))$ \textit{is cyclic}. \\
\smallskip\\
\textbf{Theorem 2.5.}($\left[ 7\right],$$\left[ 8\right]$). \textit{Let} $n$$\in$$\mathbb{N},$ $n\geq2,$ \textit{and}  $\mathbb{Q}$$\subset$$\textbf{K}$ \textit{be an extension of fields}, $[\textbf{K}:\mathbb{Q}]=n$ \textit{and} $p$
\textit{be a prime number. Let} $\mathbb{Z}_{K}$ \textit{be the ring of integers of the field} $K.$ \textit{There exist positive integers} $e_{i}$, $i=1,2,...,g,$ \textit{such that}\\
$$p\mathbb{Z}_{K}=\prod^{g}_{i=1}P^{e_{i}}_{i},$$\\
\textit{where all}$P_{i}$, $i=1,2,...,g$, \textit{are prime ideals above} $p$.\\
 \textit{The integer} $e_{i}$ \textit{is called the} \textbf{ramification index of} $p$ \textit{at} $P_{i}$. \textit{The degree} $f_{i}$  \textit{of the field extension defined by}\\
  $$f_{i}=[\mathbb{Z}_{K}/P_{i}:\mathbb{Z}/p\mathbb{Z}]$$\\
  \textit{is called} \textbf{the residual degree of p}.\\
  \smallskip\\
\textbf{Theorem 2.6.}($\left[ 7\right],$$\left[ 8\right]$).\textit{We have the following formulas:}\\
$$N(P_{i})=p^{f_{i}},\textit{and}\sum^{g}_{i=1}e_{i}f_{i}=n=[\textbf{K}:\mathbb{Q}].$$\\
\textit{In the case when} $\mathbb{Q}$$\subseteq$\textbf{K} \textit{is a Galois extension, the result is more specific}:
\textit{ the ramification indices} $e_{i}$ \textit{of} $P_{i}$ $i=1,2,...,g,$ \textit{are equal (say to e), the residual degrees} $f_{i}$ \textit{are equal as well (say to f)} \textit{and} $efg=n$. 
 \section{Results}
 \label{Results}
Now we consider the Diophantine equation\\
\begin{equation}\label{ec3.1}
   x^{4}-q^{4}=py^{r}
\end{equation}
  in the conditions
  \begin{equation}\label{ec3.2}
  p,q,r\: are\: primes,\:p\neq q\neq r\neq p, q\neq2, p\equiv3\pmod4, p\equiv1\pmod r,
  \end{equation}
   $r\equiv \pm 3 \pmod8$, $\overline{p}$ is a generator of the group $\left(U(\mathbb{Z}_{q^{r-1}}),\cdot\right)$, $\overline{q}$ is a generator of the group $\left(\mathbb{Z}_{r}^{*},\cdot\right)$, $2$ is an $r$-power residue mod $q.$\\
\smallskip\\  
We are working in the Kummer fields $\mathbb{Q}\left(\zeta;\sqrt[r]{p}\right)$ and $\mathbb{Q}\left(\zeta;\sqrt[r]{2^{r-2}p}\right)$, where $p,r$ are prime numbers, $p\equiv1$$\pmod r$.\\
\smallskip\\ 
Here are some examples of primes $p,q,r$ satsfying the conditions (3.2).\\
\smallskip\\
 \textbf{1.} $p=19$, $r=3$, $q=11$. We have: $19\equiv3$$\pmod4$, $19\equiv1$$\pmod3$, $3\equiv3$$\pmod8$, $\overline{19}$ is a generator of the group $\left(U(\mathbb{Z}_{121}),\cdot\right)$, $\overline{11}=\overline{2}$ is a generator of the group $\left(\mathbb{Z}_{3}^{*},\cdot\right)$, $2$ is a $3$-power residue mod $11$.\\
 \smallskip\\
 \textbf{2.} $p=67$, $r=3$, $q=5$. We have: $67\equiv3$$\pmod4$, $67\equiv1$$\pmod3$, $3\equiv3$$\pmod8$, $\overline{67}$ is a generator of the group $\left(U(\mathbb{Z}_{25}),\cdot\right)$, $\overline{5}=\overline{2}$ is a generator of the group $\left(\mathbb{Z}_{3}^{*},\cdot\right)$, $2$ is a $3$-power residue mod $5$.\\
 \smallskip\\
 \textbf{3.} $p=11$, $r=5$, $q=3$. We have: $11\equiv3$$\pmod4$, $11\equiv1$$\pmod5$, $5\equiv-3$$\pmod8$, $\overline{11}$ is a generator of the group $\left(U(\mathbb{Z}_{81}),\cdot\right)$, $\overline{3}$ is a generator of the group $\left(\mathbb{Z}_{5}^{*},\cdot\right)$, $2$ is a $5$-power residue mod $3$.\\
 \smallskip\\
 \textbf{4.} $p=67$, $r=11$, $q=13$. We have: $67\equiv3$$\pmod4$, $67\equiv1$$\pmod11$, $11\equiv3$$\pmod8$, $\overline{67}$ is a generator of the group $\left(U(\mathbb{Z}_{13^{10}}),\cdot\right)$, $\overline{13}=\overline{2}$ is a generator of the group $\left(\mathbb{Z}_{11}^{*},\cdot\right)$, $2$ is a $11$-power residue mod $13$.\\
 \smallskip\\ 
Let $A$ be the ring of integers of the Kummer field $\mathbb{Q}\left(\zeta;\sqrt[r]{p}\right).$
We give a general lemma about two ideals generated in $A$ by elements of the form: $y_{2}-\zeta^{m}\sqrt[r]{p}y_{1}$, with $y_{1},y_{2}$$\in$$\mathbb{Z}$ with different exponents $m\in$ $\mathbb{N}.$
\smallskip\\
\textbf{Lemma 3.1.}\textit{Let} $p$ \textit{and} $r$ \textit{be prime integers}, $p\equiv1$ (mod $r$)
\textit{and let} $\zeta$ \textit{be a primitive} $r$\textit{-th root of unity}. \textit{Let A be the ring of integers of the Kummer field} $\mathbb{Q}\left(\zeta;\sqrt[r]{p}\right),$ $y_{1}$ \textit{and} $y_{2}$ \textit{are integers such that}\:$\gcd\left(y_{1},y_{2}\right)=1$, $p$ \textit{does not divide} $y_{2}$, $m,n$$\in$$\left\{0,1,...,r-1\right\}$, $m\neq n.$ \textit{If} $y_{2}-y_{1}$ \textit{is not divisible by} $r$, 
\textit{then}\\
$$\left(y_{2}-\zeta^{m}y_{1}\sqrt[r]{p}\right)A \  \textit{and}\ \left(y_{2}-\zeta^{n}y_{1}\sqrt[r]{p}\right)A$$\\
 \textit{are coprime ideals of A}.\\
 \smallskip\\
  \textbf{Proof.}   
  We supose that $m<n$.\\
  Let $J$ be the ideal of A generated by $y_{2}-\zeta^{m}y_{1}\sqrt[r]{p}$ and $y_{2}-\zeta^{n}y_{1}\sqrt[r]{p}$.\\
 Using the Fermat's Little Theorem, we have: $y^{r}_{1}$$\equiv$$y_{1}$$\pmod r$ and 
 $y^{r}_{2}$$\equiv$$y_{2}$$\pmod r$. This implies that $y^{r}_{2}-py^{r}_{1}$$\equiv$$y_{2}-py_{1}\pmod r.$ But $p$$\equiv$$1$$\pmod r$, therefore $y^{r}_{2}-py^{r}_{1}$$\equiv$$y_{2}-y_{1}$$\pmod r.$\\
 Using the fact that $y_{2}-y_{1}$ is not divisible by $r$, it results that
  $y^{r}_{2}-py^{r}_{1}$ is not divisible by $r$, therefore gcd$(r,y^{r}_{2}-py^{r}_{1} )= 1$. This implies that there exist $h,k$$\in$ $\mathbb{Z}$ such that\\
  $$h(y^{r}_{2}-py^{r}_{1})+ kr=1. $$\\
  In the ring $A$, we have\\
 $$y^{r}_{2}-py^{r}_{1}=(y_{2}-y_{1}\sqrt[r]{p})(y_{2}-y_{1}\zeta\sqrt[r]{p})...
  (y_{2}-y_{1}\zeta^{r-1}\sqrt[r]{p}).$$\\
  This implies that:\\
  $$y^{r}_{2}-py^{r}_{1}\in J. $$\\
 Since\\ $$(y_{2}-y_{1}\zeta^{m}\sqrt[r]{p})-(y_{2}-y_{1}\zeta^{n}\sqrt[r]{p})=\zeta^{m}y_{1}\sqrt[r]{p}(\zeta^{n-m}-1)=\zeta^{m}y_{1}\sqrt[r]{p}u_{n-m}(\zeta-1),$$ where $u_{n-m}$, $\zeta^{m}$$\in$$U(\textbf{Z}[\zeta])$$\subset$$U(A),$ we obtain that $y_{1}\sqrt[r]{p}(\zeta-1)$$\in$J. But $\sqrt[r]{p^{r-1}}$$\in$A, therefore $py_{1}(\zeta-1)$$\in$J. Now,
  $y_{2}-y_{1}\zeta^{m}\sqrt[r]{p}$$\in$J and $\zeta^{n-m}$$\in$A. These imply that $y_{2}\zeta^{n-m}-\zeta^{n}y_{1}\sqrt[r]{p}$$\in$J. Using the fact that $y_{2}-y_{1}\zeta^{n}\sqrt[r]{p}$$\in$J, it results that $y_{2}(\xi^{n-m}-1)$$\in$J. But $\zeta^{n-m}-1=u_{n-m}(\zeta-1)$, with $u_{n-m}$$\in$U(A). We obtain that
  $y_{2}(\zeta-1)$$\in$J.\\ 
  From the hypothesis, we know that gcd$(y_{1},y_{2})=1$ and $y_{2}$ is not divisible with $p$. This implies that $\gcd(py_{1},y_{2})=1$, therefore there exist 
  $h_{1},$ $h_{2}$$\in$$\mathbb{Z}$ such that $py_{1}h_{1}+y_{2}h_{2}=1$.
   Multiplying the last equality by $\zeta-1$ and using the previous relations, we obtain that $\zeta-1$$\in$J.\\
   Since
   $r=u(1-\zeta)^{r-1}$, where $u$$\in$$U(\mathbb{Z}[\zeta])$$\subset$$U(A)$, we have $r$$\in$J.\\  
   From this relation and the previous ones, it results that $1$$\in$J, therefore 
   $$\left(y_{2}-\zeta^{m}y_{1}\sqrt[r]{p}\right)A \quad\makebox{and} \quad \left(y_{2}-\zeta^{n}y_{1}\sqrt[r]{p}\right)A$$\\
 are coprime ideals of A.\\
 \smallskip\\
 Analogously we may prove:
\smallskip\\
\textbf{Lemma 3.2.}\textit{Let} $p$ \textit{and} $r$ \textit{be prime integers}, $p\equiv1$ $(mod r)$
\textit{and let} $\zeta$ \textit{be a primitive} $r$\textit{-th root of unity}. \textit{Let A be the ring of integers of the Kummer field} $\mathbb{Q}\left(\zeta;\sqrt[r]{2^{r-2}p}\right),$ $y_{1}$ \textit{and} $y_{2}$ \textit{integers such that} $\gcd\left(y_{1},y_{2}\right)=1$, $p$ \textit{does not divide} $2y_{2}$, $m,n$$\in$$\left\{0,1,...,r-1\right\}$, $m\neq n.$ \textit{If} $y_{2}-2^{r-2}py_{1}$ \textit{is not divisible by} $r$,  
\textit{then} \\
$$\left(y_{2}-\zeta^{m}y_{1}\sqrt[r]{2^{r-2}p}\right)A \  \textit{and}\ \left(y_{2}-\zeta^{n}y_{1}\sqrt[r]{2^{r-2}p}\right)A$$\\
 \textit{are coprime ideals of A}.
  \smallskip\\
  We may consider now our equation $x^{4}-q^{4}=py^{r}$.\\
  \smallskip\\
  \textbf{Proof of the Main theorem.}We reason by reduction to absurd.Let $\left(x,y\right)$$\in$$\mathbb{Z}^{2},$ $\gcd(x,y)=1$ be a solution to the equation \:(3.1), with $p,q,r$ satisfying conditions (3.2).\\
  Suppose, by way of contradiction that $p$ does not divide $y.$\\
  We consider two cases: either \textit{x} is odd or \textit{x} is even.\\
  \textbf{Case 1}. $\textbf{x}$ \textbf{is an odd number}\\
  Since \textit{q} is a prime number, $q\geq3$, we get $x^{2},q^{2}\equiv1$ $\pmod4,$ therefore $x^{2}-q^{2}\equiv0$ $\pmod4$, $x^{2}+q^{2}\equiv2$ $\pmod4$.\\
 We denote $d=$$\gcd\left(x^{2}-q^{2},x^{2}+q^{2}\right)$. Then $d|2x^{2}$ and  $d|2q^{2}$. But $\gcd\left(x,y\right)=1$ implies that $x$ is not divisible by $q$. Therefore $d=2$. We get either\\
 \begin{center}
$x^{2}-q^{2}=2^{r-1}py^{r}_{1}$ and $x^{2}+q^{2}=2y^{r}_{2},$
\end{center}
where $y_{1},y_{2}$$\in$$\mathbb{Z}$, $2y_{1}y_{2}=y$, $y_{2}$ is an odd number, $\gcd\left(y_{1},y_{2}\right)=1,$\\ 
or \\
\begin{center}
$x^{2}-q^{2}=2^{r-1}y^{r}_{1}$ and $x^{2}+q^{2}=2py^{r}_{2},$
\end{center}
where $y_{1},y_{2}$$\in$$\mathbb{Z}$, $2y_{1}y_{2}=y$, $y_{2}$ is an odd number, $\gcd\left(y_{1},y_{2}\right)=1.$\\
In the last case, we obtain that $p|\left(x^{2}+q^{2}\right)$, in contradiction with the fact that $p\equiv3$$\pmod4$.
It remains to study the case 
\begin{center}
$x^{2}-q^{2}=2^{r-1}py^{r}_{1}$ and $x^{2}+q^{2}=2y^{r}_{2}$.
\end{center}
By subtracting the two equations, we obtain $q^{2}=y^{r}_{2}-2^{r-2}py^{r}_{1}$.\\
We consider the Kummer field $\mathbb{Q}\left(\zeta,\sqrt[r]{2^{r-2}p}\right)$, where $\zeta$ is a primitive $r-$ th root of unity.\\
 The last equality becomes:\\ $$q^{2}=\left(y_{2}-y_{1}\sqrt[r]{2^{r-2}p}\right)\left(y_{2}-y_{1}\zeta\sqrt[r]{2^{r-2}p}\right)...\left(y_{2}-y_{1}\zeta^{r-1}\sqrt[r]{2^{r-2}p}\right)$$\\
 But $<\overline{q}>=$$(\mathbb{Z}^{*}_{r}, \cdot)$ and, by applying Proposition 2.1, we obtain that $q\mathbb{Z}\left[\zeta\right]$ is a prime ideal in the ring $\mathbb{Z}\left[\zeta\right]$.\\ 
   We try to decompose the ideal $(q)$ in the ring $A.$ We have:\\  $$\left\{\frac{2^{r-2}p}{\left(q\right)}\right\}=\left\{\frac{2}{\left(q\right)}\right\}^{r-2}\left\{\frac{p}{\left(q\right)}\right\}.$$\\
 Since $2$ is an $r$-power residue mod $q$, then there is $\alpha$$\in$$\mathbb{Z}\left[\zeta\right]$ such that $\alpha^{r}$$\equiv$$2$ $\pmod q$, therefore $\left\{\frac{2}{\left(q\right)}\right\}=1$.\\
 We obtain that $\left\{\frac{2^{r-2}p}{\left(q\right)}\right\}=\left\{\frac{p}{\left(q\right)}\right\}$
   and we get:\\ $$\left\{\frac{p}{\left(q\right)}\right\}=\zeta^{c}\equiv\\p^{\frac{N\left((q)\right)-1}{r}}\pmod q.$$\\
   We next calculate $N\left((q)\right)$.\\
   Since $q\mathbb{Z}\left[\zeta\right]$ is a prime ideal in the ring $\mathbb{Z}\left[\zeta\right]$ it results that $e=1$, $g=1$. But $efg=[\mathbb{Q}(\zeta):\mathbb{Q}]=r-1$, therefore $f=r-1$. Using Theorem 2.6, we obtain $N\left((q)\right)=q^{r-1}$
   and\\ $$\left\{\frac{p}{\left(q\right)}\right\}=\zeta^{c}\equiv\\p^{\frac{q^{r-1}-1}{r}}\pmod q.$$\\
   If $\left\{\frac{p}{\left(q\right)}\right\}=1$, it results $p^{\frac{q^{r-1}-1}{r}}\equiv1$$\pmod q$. But\\
   $<\overline{p}>=(U(Z_{q^{r-1}});\cdot)$ and $\left|U(Z_{q^{r-1}})\right|=q^{r-1}-q^{r-2}$, hence $q^{r-1}-q^{r-2}|\frac{q^{r-1}-1}{r}$. This implies that there is $j\in\mathbb{N}^{*}$ such that $\frac{q^{r-1}-1}{r}=j(q^{r-1}-q^{r-2})$. Therefore $q^{r-2}+q^{r-3}+...+q+1=jrq^{r-2}$. This equality is impossible, because, for $q\in\mathbb{N}^{*}$, $q\geq3$, we have:
   $q^{r-2}+q^{r-3}+...+q+1\leq(r-1)q^{r-2}<jrq^{r-2}$.\\
   We obtain that 
$\left\{\frac{p}{\left(q\right)}\right\}=\zeta^{c}\neq1$, therefore $qA\in$Spec(A).
  Passing to ideals in the expresion of $q^2$, we get:\\ $$\left(y_{2}-y_{1}\sqrt[r]{2^{r-2}p}\right)A\left(y_{2}-y_{1}\zeta\sqrt[r]{2^{r-2}p}\right)A...\left(y_{2}-y_{1}\zeta^{r-1}\sqrt[r]{2^{r-2}p}\right)A=\left(qA\right)^{2}\ \ \ $$\\
   and,according to Lemma 3.2, this equality is impossible.\\
   \textbf{Case 2.} $\textbf{x}$ \textbf{ is an even number.}\\   
In this case, $x^{2}-q^{2}$ and $x^{2}+q^{2}$ are odd numbers.\\
We prove that $\gcd\left(x^{2}-q^{2},x^{2}+q^{2}\right)=1$. We suppose that there exists an odd prime number $d$ such that $d|\left(x^{2}-q^{2}\right)$ and $d|\left(x^{2}+q^{2}\right)$. Hence $d|x$ and $d|q$. Using the hypothesis, we obtain that $d|y$, in contradiction with the fact $\gcd\left(x,y\right)=1$.  Therefore gcd$\left(x^{2}-q^{2},x^{2}+q^{2}\right)=1$.
The equation(3.1) implies either\\
\begin{center}
$x^{2}-q^{2}=py^{r}_{1},$ $x^{2}+q^{2}=y^{r}_{2}$,
with $y_{1},y_{2}$$\in$$\mathbb{Z}$, $y_{1}y_{2}=y,$ $\gcd\left(y_{1},y_{2}\right)=1$
\end{center}
 or\\
\begin{center}
$x^{2}-q^{2}=y^{r}_{1},$ $x^{2}+q^{2}=py^{r}_{2}$,
with $y_{1},y_{2}$$\in$$\mathbb{Z}$, $y_{1}y_{2}=y,$ $\gcd\left(y_{1},y_{2}\right)=1.$\\
\end{center}
In the last case, we obtain that $p|\left(x^{2}+q^{2}\right)$, in contradiction with the fact that $p\equiv3$$\pmod4$.
It remains the case 
\begin{center}
$x^{2}-q^{2}=py^{r}_{1}$ and $x^{2}+q^{2}=y^{r}_{2}.$
\end{center}
Subtracting the two equations, we get $2q^{2}=y^{r}_{2}-py^{r}_{1}$.\\
Let $\mathbb{Q}\left(\zeta,\sqrt[r]{p}\right)$ be a Kummer field, where $\zeta$ is a primitive $r-$ th root of unity, and A the ring of integers of $\mathbb{Q}\left(\zeta,\sqrt[r]{p}\right).$ In A, the last equality becomes
$$\left(y_{2}-y_{1}\sqrt[r]{p}\right)\left(y_{2}-y_{1}\zeta\sqrt[r]{p}\right)...\left(y_{2}-y_{1}\zeta^{r-1}\sqrt[r]{p}\right)=2q^{2}.$$ 
   We prove that $2$ is a prime element in the ring $\mathbb{Z}[\zeta]$.\\
   According to Proposition 2.1, $2\mathbb{Z}\left[\zeta\right]=P_{1}P_{2}...P_{s}$, where  $s=\frac{\varphi(r)}{ord_{\textbf{Z}^{*}_{r}}(\overline{2})}=\frac{r-1}{ord_{\textbf{Z}^{*}_{r}}(\overline{2})}$; but
   $r\equiv3$ or $5$$\pmod8$ implies that $2$ is not a quadratic residue modulo $r$. Applying the Euler's Criterion, we obtain that $2^{\frac{r-1}{2}}\equiv-1$ $\pmod r$. Hence $ord_{\mathbb{Z}^{*}_{r}}(\overline{2})=r-1.$\\
   We get that $2\mathbb{Z}\left[\zeta\right]$$\in$Spec$(\mathbb{Z}\left[\zeta\right])$\\
As $p$ is a prime number which satisfies $p$$\equiv$$3$$\pmod4$, we have $\left\{\frac{p}{\left(2\right)}\right\}=1$. Using Theorem 2.2, we obtain that $2A=P^{'}_{1}P^{'}_{2}...P^{'}_{r}$, where $P^{'}_{1},P^{'}_{2},...,P^{'}_{r}$ are prime ideals in the ring $A$. \\
As in the case 1, $qA$$\in$Spec($A$).
By considering ideals in the relation in the expression of $2q^{2}$, we obtain\\ 
\begin{equation}\label{ec3} \left(y_{2}-y_{1}\sqrt[r]{p}\right)A\left(y_{2}-y_{1}\zeta\sqrt[r]{p}\right)A...\left(y_{2}-y_{1}\zeta^{r-1}\sqrt[r]{p}\right)A=P^{'}_{1}P^{'}_{2}...P^{'}_{r}(qA)^{2}.
\end{equation}
Let G be the Galois group of the Kummer field $\mathbb{Q}(\zeta,\sqrt[r]{p})$ over $\mathbb{Q}(\zeta)$. According to Proposition 2.4, G is a cyclic group with $\sigma$ as a generator, $\sigma:\mathbb{Q}(\zeta,\sqrt[r]{p})\rightarrow\mathbb{Q}(\zeta,\sqrt[r]{p})$, 
$\sigma(\sqrt[r]{p})=\zeta\sqrt[r]{p}$.
  We come back to the equality (3.3) and we consider three cases.
\smallskip\\
\textbf{Subcase (i):} If there exists $k$$\in$$\left\{1,2,...,r\right\}$ such that $\left(y_{2}-y_{1}\sqrt[r]{p}\right)A=P^{'}_{k}$$\in$Spec(A), we use Proposition 2.3, and we obtain that $\sigma\left(\left(y_{2}-y_{1}\sqrt[r]{p}\right)A\right)=\left(y_{2}-y_{1}\zeta\sqrt[r]{p}\right)A$$\in$Spec(A) and $\sigma^{2}\left(\left(y_{2}-y_{1}\sqrt[r]{p}\right)A\right)=\left(y_{2}-y_{1}\zeta^{2}\sqrt[r]{p}\right)A$$\in$Spec(A),...,\\ $\sigma^{r-1}\left(\left(y_{2}-y_{1}\sqrt[r]{p}\right)A\right)=\left(y_{2}-y_{1}\zeta^{r-1}\sqrt[r]{p}\right)A$$\in$Spec(A),
therefore the equality (3.3) is impossible.
\smallskip\\
\textbf{Subcase (ii):} If $\left(y_{2}-y_{1}\sqrt[r]{p}\right)A=P^{2}_{11}$, where $P_{11}$$\in$ Spec(A), using Proposition 2.3, we obtain $\left(y_{2}-y_{1}\zeta\sqrt[r]{p}\right)A=\sigma(\left(y_{2}-y_{1}\sqrt[r]{p}\right)A )=P^{2}_{12}$, where $P_{12}$$\in$ Spec(A), and so on, up to 
$\left(y_{2}-y_{1}\zeta^{r-1}\sqrt[r]{p}\right)A=\sigma^{r-1}(\left(y_{2}-y_{1}\sqrt[r]{p}\right)A )=P^{2}_{1r}$, where $P_{1r}$$\in$ Spec(A).\\
The equality (3.3) becomes\\
$$P^{2}_{11}P^{2}_{12}.....P^{2}_{1r}=P^{'}_{1}P^{'}_{2}...P^{'}_{r}(qA)^{2},$$\\
where $P_{1i}$$\in$ Spec(A), $(\forall)i=\overline{1,r}$ $P^{'}_{i}$$\in$ Spec(A), $i=1,2,...,r$, $qA$$\in$ Spec(A). This equality is impossible, since all intervening ideals are prime in $A$ and $r>2.$\\
\smallskip\\
\textbf{Subcase (iii):} If 
$\left(y_{2}-y_{1}\sqrt[r]{p}\right)A=P_{21}P^{'}_{21},$ where  $P_{21}$,$P^{'}_{21}$$\in$Spec(A), $P_{21}$$\neq$$P^{'}_{21},$ then\\ $\sigma\left(\left(y_{2}-y_{1}\sqrt[r]{p}\right)A\right)=\left(y_{2}-y_{1}\zeta\sqrt[r]{p}\right)A=P_{22}P^{'}_{22},$ where  $P_{21}$,$P^{'}_{21}$ are distinct prime ideals in the ring A and so on, up to we have in the end\\
$$\sigma^{r-1}\left(\left(y_{2}-y_{1}\sqrt[r]{p}\right)A\right)=\left(y_{2}-y_{1}\zeta^{r-1}\sqrt[r]{p}\right)A=P_{2r}P^{'}_{2r},$$ 
where $P_{2r},P^{'}_{2r}$\;are\; prime\; ideals \;in \;A. \\
Therefore relation (3.3) is equivalent to\\ $$P_{21}P^{'}_{21}P_{22}P^{'}_{22}.....P_{2r}P^{'}_{2r}=P^{'}_{1}P^{'}_{2}...P^{'}_{r}(qA)^{2},  \\$$\\
 where $P_{21}$,$P^{'}_{21}$,$P_{22}$,$P^{'}_{22}$,.....,$P_{2r}$,$P^{'}_{2r}$,$P^{'}_{1}$,$P^{'}_{2}$...,$P^{'}_{r}$,$qA$ \;are\; prime\; ideals \;in \;A.\\ This equality does not hold.\\
  From subcases (i), (ii), (iii), it results that the equality (3.3) is impossible.\\
  We get that the supposition made is false, so $p$ must divide $y.$\\
  \smallskip\\
  \textbf{Concluding remarks}\\
  Finally we give some examples where at least one of the conditions (3.2) is not ful filled and Main Theorem is false. \\
  \smallskip\\
  \textbf{1.} $p=17$, $r=5$, $q=3$. We have: $5\equiv-3$$\pmod8$, $\overline{3}$ is a generator of the group $\left(\mathbb{Z}_{5}^{*},\cdot\right),$ $2$ is a $5$-power residue mod $3,$ but $17$ is not congruent with $3$$\pmod4$, $17$ is not congruent with $1$ $\pmod3$,$\overline{17}$ is not a generator of the group $\left(U(\mathbb{Z}_{81}),\cdot\right)$.\\
  In this situation we observe that $(x,y)=(5,2)$ is a solution of the equation
  $x^{4}-3^{4}=17y^{5}$ but $p$ does not divide $y.$\\
  \smallskip\\
  \textbf{2.} $p=65537$, $r=11$, $q=257$. We have: $11\equiv3$$\pmod8$, $\overline{257}=\overline{4}$ is not a generator of the group $\left(\mathbb{Z}_{11}^{*},\cdot\right)$, $2$ is a $11$-power residue mod $257,$ $\overline{65537}$ is not a generator of the group $\left(U(\mathbb{Z}_{257^{10}}),\cdot\right)$, $65537$ is not congruent with $3$$\pmod4$, $65537$ is not congruent with $1$$\pmod r$.\\
  In this situation we observe that $(x,y)=(255,-2)$ is a solution of the equation
  $x^{4}-257^{4}=65537y^{11}$ but $p$ does not divide $y.$\\
  \smallskip\\
  \textbf{3.} $p=257$, $r=7$, $q=17$. We have: $7$ is not congruent with $\pm3$ $\pmod8$, $\overline{17}=\overline{3}$ is a generator of the group $\left(\mathbb{Z}_{7}^{*},\cdot\right)$, $2$ is a $7$-power residue mod $17,$ $\overline{257}$ is not a generator of the group $\left(U(\mathbb{Z}_{17^{6}}),\cdot\right)$, $257$ is not congruent with $3$ $\pmod4$, $257$ is not congruent with $1$$\pmod7$.\\
  In this situation we observe that $(x,y)=(15,-2)$ is a solution of the equation
  $x^{4}-17^{4}=257y^{7}$ but $p$ does not divide $y.$\\
  \smallskip\\
  In the last example we have $r\equiv-1$ (mod $8$).This situation is not covered by Main Theorem. The particular case $r=7$ has been subject of investigation in [14]. In near future we intend to analyse equation 3.1 also in the case $r\equiv\pm1 \pmod8.$
  \smallskip\\
  Another remark is that there is a connection between the last three examples, more exactly, if $F_{n}$ is the $n$-term of Fermat sequence ($F_{n}=2^{2^{n}}+1$), in example \textbf{1.} we have $p=17=F_{2};$ in example \textbf{3.} we have $q=17=F_{2},$ $p=257=F_{3};$ in example \textbf{2.} we have $q=257=F_{3},$ $p=65537=F_{4}.$\\
  \smallskip\\
  In the future we will study if Lemma 3.1 and Lemma 3.2 are valid not only for $p\equiv1$(mod $r$) but also in more general conditions.\\
  Another research theme is to investigate to what extent the conditions imposed on $p,$ $q,$ $r$ can be relaxed, so that analogous results could be obtained in more general conditions.\\       
  \bigskip\\
  \textbf{References}
  \smallskip\\
  $\left[1\right]$ H. Cohen, \textit{A Course in Computational Algebraic Number Theory}, Springer-Verlag, 1993.\\
  $\left[2\right]$ D. Cox, \textit{Primes of the Form} $x^{2}+ny^{2}$, A. Wiley -Interscience Publication, New York, 1989.\\
  $\left[3\right]$ H. Darmon, \textit{The equation} $x^{4}-y^{4}=z^{p},$ C.R. Math. Rep. Acad. Sci. Canada. XV No. 6 (1993), 286-290.\\
  $\left[4\right]$ K. Gy\H ory, \textit{Uber die diophantische Gleichung} $x^{p}+y^{p}=cz^{p},$ Publ. Math. Debrecen, 13 (1966), 301-306.\\
  $\left[5\right]$ K. Gy\H ory, \textit{On the diophantine equation} $x^{p}+y^{p}=cz^{p},$ Mat. Lapok, 18 (1967), 93-96.\\
  $\left[6\right]$ K. Gy\H ory, A. Peth\H o and V. T. S$\acute{o}$s, \textit{Number Theory, Diophantine, Computational and Algebraic Aspects,} Walter de Gruyter, Berlin-New York, 1998.\\
  $\left[7\right]$ D.Hilbert, \textit{The theory of algebraic number fields,} Springer-Verlag,1998.\\
  $\left[8\right]$ F. Lemmermeyer, \textit{Reciprocity Laws}, Springer-Verlag, 2000.\\
  $\left[9\right]$ F. Luca, A.Togbe, \textit{On the Diophantine Equation} $x^{4}-q^{4}=py^{3},$ Accepted.\\
  $\left[10\right]$ B.Powell, \textit{Sur l'$\acute{e}$quation Diophantienne}  $x^{4}$$\pm$$y^{4}=z^{p},$ Bull. Sc. Math., 107(1983), 219-223.\\
  $\left[11\right]$ P.Ribenboim, \textit{Classical Theory of Algebraic Numbers}, Universitext,Springer, 2001.\\
  $\left[12\right]$ D. Savin, \textit{On the Diophantine Equation} $x^{4}-q^{4}=py^{3}$, An. St. Univ.Ovidius, Ser. Mat., \textbf{12}(2004), fasc.1, p. 81 - 90.\\
  $\left[13\right]$ D. Savin, \textit{On the Diophantine Equation} $x^{4}-q^{4}=py^{5}$, Italian Journal of Pure and Applied Mathematics (paper accepted for publication).\\
  $\left[14\right]$ D. Savin, A. Barbulescu \textit{The Diophantine Equation} $x^{4}-q^{4}=py^{7}$ \textit{in Special Conditions}, Journal Automation Computers Applied Mathematics, vol. \textbf{15}(2006), no.2, 295-300.\\
  \bigskip\\
{Faculty of Mathematics and Computer Science ,\\
  Department of Mathematics\\
"Ovidius" University of Constanta}\\
Bd. Mamaia 124, Constanta, 900527\\
Romania\\
{e-mail: Savin.Diana@univ-ovidius.ro\\
\ \ \ \ \ dianet72@yahoo.com} 
 \end{document}